\newmathalphabet*{\E}{eus}{m}{n}

\renewcommand{\P}{{\Bbb P}}

\newcommand{\F}{{\E F}}

\renewcommand{\O}{{\E O}}
\newcommand{\Z}{{\Bbb Z}}
\newcommand{\R}{{\Bbb R}}
\newcommand{\N}{{\Bbb N}}

\newcommand{\Spec}{\operatorname{Spec}}

\newcommand{\codim}{\operatorname{codim}}

\newcommand{\mod}{\operatorname{mod}}
\newcommand{\Hom}{\operatorname{Hom}}

\renewcommand{\H}{\operatorname{H}}

\newcommand{\Pic}{\operatorname{Pic}}

\newtheorem{thm}{Theorem}
\newtheorem{lemma}{Lemma}
\newtheorem{prop}{Proposition}
\newtheorem{cor}{Corollary}

\newtheorem{remark}{Remark}

\newenvironment{pf}{{\bf Proof}}{$\Box$}

\newcommand{\binom}[2]{{\left(\begin{array}{c} #1\\#2 \end{array}\right)}}

\documentstyle[amsfonts,amscd,12pt]{amsart}

\title{Frobenius morphisms over $\Z/p^2$ and Bott vanishing}

\date{\today}

\author{Anders Buch}
\address[Anders Buch]
{Matematisk Institut\\
Aarhus Universitet\\
Ny Munkegade\\
DK-8000 \AA rhus C\\
Denmark}
\email[Anders Buch]{abuch@@mi.aau.dk}
\author{Jesper Funch Thomsen}
\address[Jesper Funch Thomsen]
{Matematisk Institut\\
Aarhus Universitet\\
Ny Munkegade\\
DK-8000 \AA rhus C\\
Denmark}
\email[Jesper Funch Thomsen]{funch@@mi.aau.dk}
\author{Niels Lauritzen}
\address[Niels Lauritzen]
{Matematisk Institut\\
Aarhus Universitet\\
Ny Munkegade\\
DK-8000 \AA rhus C\\
Denmark}
\email[Niels Lauritzen]{niels@@mi.aau.dk}
\author{Vikram Mehta}
\address[Vikram Mehta]
{School of Mathematics\\
Tata Institute of Fundamental Research\\
Homi Bhabha Road\\
Bombay}
\email[Vikram Mehta]{vikram@@tifrvax.tifr.res.in}



\begin{document}

\maketitle

In \cite{Bott} Bott proved the vanishing theorem
$$
\H^\bullet(\P^n, \Omega_{\P^n}^p\otimes \O(k))=0
$$
when $0<k\leq p$. We will say that a smooth projective
variety $X$ has Bott vanishing if for every ample line
bundle $L$, $i>0$ and $j$
$$
\H^i(X, \Omega_X^j\otimes L)=0
$$

The purpose of this paper is to show that Bott vanishing is a simple
consequence of a very specific condition on the Frobenius morphism in prime
characteristic $p$.
Recall that the absolute Frobenius morphism
$F:X\rightarrow X$ on $X$, where $X$ is a variety over $\Z/p$
is the identity on point
spaces and the $p$-th power map locally on functions. Assume that
there is a flat scheme $X^{(2)}$ over $\Z/p^2$, such that $X\cong X^{(2)}
\times_{\Z/p^2}
\Z/p$. The condition on $F$ is that there should be a morphism
$F^{(2)}:X^{(2)}\rightarrow X^{(2)}$ which gives $F$ by reduction mod $p$.
In this case we will say that the Frobenius morphism lifts to $\Z/p^2$. It is
known that a lift of the Frobenius morphism to $\Z/p^2$ leads to
a quasi-isomorphism
$$
\sigma:\bigoplus_{0\leq i}\Omega^i_X[-i] @>>> F_* \Omega^\bullet_X
$$
where the complex on the left has zero differentials and $\Omega^\bullet_X$
denotes the de Rham complex of $X$ (\cite{DI}, Remarques 2.2(ii)).
Using duality we prove that $\sigma$ is in fact a split
quasi-isomorphism.

In general it is very difficult to decide when Frobenius lifts to $\Z/p^2$.
However for varieties which are glued together by monomial automorphisms
it is easy. This is the case for toric varieties, where we show
that the Frobenius morphism lifts to $\Z/p^2$. This gives natural
characteristic $p$ proofs and explanations of the Bott vanishing theorem for
(singular and smooth) toric varieties and the degeneration of the Danilov
spectral sequence (\cite{Dan}, Theorem 7.5.2, Theorem 12.5).

Paranjape and Srinivas have
proved using complex algebraic geometry that if Frobenius for a generalized
flag variety $X$ lifts to the $p$-adic numbers $\hat{\Z}_p=
\underleftarrow{\operatorname{lim}}_n \Z/p^n$, then $X$ is a product of
projective spaces \cite{PaSri}.
In the last part of this paper we generalize this
result by showing that Frobenius for a large class of generalized flag
varieties admits no lift to $\Z/p^2$. This is done using a lemma
on fibrations linking non-lifting of Frobenius to Bott
non-vanishing cohomology groups for flag varieties of Hermitian
symmetric type over the complex numbers. These cohomology groups
have been studied thoroughly by M.-H.~Sato and D.~Snow.
\vskip 0.8truecm
This paper grew out of a seminar on characteristic $p$ methods
at Aarhus University focusing on Deligne and Illusie's algebraic
proof \cite{DI} of Kodaira's vanishing theorem. We are grateful to the
patient and inspiring participants: H.~H.~Andersen, M.~B\"okstedt,
J.~A.~Geertsen, S.~H.~Hansen, J.~C.~Jantzen and J.~Kock.
\section{Preliminaries}

Let $k$ be a perfect field of characteristic $p>0$ and
$X$ a smooth $k$-variety of dimension $n$. By $\Omega_X$ we denote
the sheaf of $k$-differentials on $X$ and $\Omega^j_X=\wedge^j \Omega_X$.
The (absolute) Frobenius morphism $F:X\rightarrow X$ is the morphism
on $X$, which is the identity on the level of points and given
by $F^\#(f)=f^p: \O_X(U)\rightarrow F_*\O_X(U)$ on the level of
functions. If $\F$ is an $\O_X$-module, then $F_* \F=\F$ as sheaves
of abelian groups, but the $\O_X$-module structure is changed according
to the homomorphism $\O_X\rightarrow F_*\O_X$.

\subsection{The Cartier operator}

The universal derivation $d:\O_X \rightarrow \Omega_X$ gives rise
to a family of $k$-homomorphisms $d^j: \Omega^j_X\rightarrow \Omega^{j+1}_X$
making $\Omega^\bullet_X$ into a complex of $k$-modules which is
called the de Rham complex of $X$. By applying
$F_*$ to the de Rham complex, we obtain a complex $F_*\Omega^\bullet_X$ of
$\O_X$-modules.
Let $B^i_X\subseteq Z^i_X\subseteq F_*\Omega^i_X$
denote the coboundaries and cocycles in degree $i$. There is the following
very nice description of the cohomology of $F_*\Omega^\bullet_X$ due
to Cartier.

\begin{thm}
There
is a uniquely determined graded $\O_X$-algebra isomorphism
$$
C^{-1}:\Omega_X^\bullet\rightarrow \cal H^\bullet(F_* \Omega^\bullet_X)
$$
which in degree $1$ is given locally as
$$
C^{-1}(da)= a^{p-1} da
$$
\end{thm}
\begin{pf}
\cite{Katz}, Theorem 7.2.
\end{pf}

With some abuse of notation, we let $C$ denote the natural homomorphism
$Z^i_X\rightarrow\Omega^i_X$,
which after reduction modulo $B^i_X$ gives the inverse isomorphism
to $C^{-1}$. The isomorphism $\bar{C}:Z^i_X/B^i_X\rightarrow \Omega^i_X$
is called the Cartier operator.

\section{Liftings of Frobenius to $W_2(k)$}
\label{flift}
There is a very interesting connection between the Cartier operator
and liftings of the Frobenius morphism to flat schemes of characteristic
$p^2$. This beautiful observation was first made by Mazur in \cite{Maz}. We
go on to explore this next.

\subsection{Witt vectors of length two}

The Witt vectors $W_2(k)$ (\cite{MumCu}, Lecture 26) of
length $2$ over $k$ can be interpreted as the set
$k\times k$, where
multiplication and addition for $a=(a_0, a_1)$ and $b=(b_0, b_1)$ in
$W_2(k)$ are defined by
$$
a\, b=(a_0\, b_0, a_0^p b_1+ b_0^p a_1)
$$
and
$$
a+b=(a_0+b_0, a_1+b_1+\sum_{j=1}^{p-1} p^{-1}\binom{p}{j} a_0^j\, b_0^{p-j})
$$
In the case $k=\Z/p$, one can prove that $W_2(k)\cong \Z/p^2$.
The projection on the first coordinate $W_2(k)\rightarrow k$ corresponds
to the reduction $W_2(k)\rightarrow W_2(k)/p\cong k$ modulo $p$.
The ring homomorphism $F^{(2)}:W_2(k)\rightarrow W_2(k)$ given by
$F^{(2)}(a_0, a_1)=(a_0^p, a_1^p)$ reduces to
the Frobenius homomorphism $F$ on $k$ modulo $p$.

\subsection{Splittings of the de Rham complex}

The previous section shows that there is a canonical morphism
$\Spec k\rightarrow \Spec W_2(k)$. Assume that there is a flat
scheme $X^{(2)}$ over $\Spec W_2(k)$ such that
\begin{equation}
\label{modp}
X\cong X^{(2)}\times_{\Spec W_2(k)}\Spec k
\end{equation}
We shall say that the Frobenius morphism $F$ lifts to $W_2(k)$ if
there exists a morphism $F^{(2)}:X^{(2)}\rightarrow X^{(2)}$ covering
the Frobenius homomorphism $F^{(2)}$ on $W_2(k)$, which
reduces to $F$ via the isomorphism (\ref{modp}). When we use the
statement that Frobenius lifts to $W_2(k)$ we will always implicitly
assume the existence of the flat lift $X^{(2)}$.

\begin{thm}
\label{split}
If the Frobenius morphism on $X$ lifts to $W_2(k)$
then there is a split quasi-isomorphism
$$
0 @>>> \bigoplus_{0\leq i}\Omega^i_X[-i] @>\sigma>> F_* \Omega^\bullet_X
$$
\end{thm}
\begin{pf}
For
an affine open subset $\Spec A^{(2)}\subseteq X^{(2)}$ there
is a ring homomorphism $F^{(2)}: A^{(2)}\rightarrow A^{(2)}$
such that
$$
F^{(2)}(b)=b^p + p\cdot \varphi(b)
$$
where $\varphi: A^{(2)}\rightarrow A=A^{(2)}/ p A^{(2)}$ is some
function and
$p\, \cdot: A\rightarrow A^{(2)}$ is the $A^{(2)}$-homomorphism
derived from
tensoring the short exact sequence of $W_2(k)$-modules
$$
\CD
0 @>>> p\,W_2(k)@>>> W_2(k) @>p\,\cdot>> p\,W_2(k) @>>> 0
\endCD
$$
with the flat $W_2(k)$ module $A^{(2)}$
identifying $A\cong A^{(2)}/p A^{(2)}$ with $p\, A^{(2)}$.
We get the following properties
of $\varphi$:
\begin{align*}
\varphi(a+b)&=\varphi(a)+\varphi(b)-\sum_{j=1}^{p-1} p^{-1} \binom{p}{j}
\bar{a}^j \bar{b}^{p-j} \\
\varphi(a\, b)&=\bar{a}^p \varphi(b)+ \bar{b}^p \varphi(a)
\end{align*}
where $\bar{\cdot}$ means reduction $\mod p$.
Now it follows that
$$
a\mapsto a^{p-1} da+ d \varphi(\tilde{a})
$$
where $\tilde{a}$ is any lift of $a$,
is a well defined derivation $\delta:A\rightarrow Z^1_{\Spec A}
\subset F_*\Omega^1_{\Spec A}$,
which gives a homomorphism $\varphi:\Omega_{\Spec A}^1\rightarrow Z^1_{\Spec A}
\subset F_* \Omega^1_{\Spec A}$. This homomorphism can be extended via the
algebra structure to give an $A$-algebra homomorphism
$\sigma: \oplus_i \Omega_{\Spec A}^i \rightarrow Z^\bullet_{\Spec A}\subseteq
F_*\Omega_{\Spec A}^\bullet$,
which composed with the
canonical homomorphism $Z^\bullet_{\Spec A}\rightarrow
\cal H^\bullet(F_* \Omega^\bullet_{\Spec A})$ gives the inverse
Cartier operator.
Since an affine open covering $\{\Spec A^{(2)}\}$ of $X^{(2)}$ gives rise
to an affine open covering $\{\Spec A^{(2)}/p A^{(2)}\}$ of $X$,
we have proved that $\sigma$ is a quasi-isomorphism of complexes
inducing the inverse Cartier operator on cohomology.

Now we give a splitting
homomorphism of $\sigma_i:\Omega^i_X\rightarrow F_*\Omega^i_X$. Notice
that $\sigma_0:\O_X\rightarrow F_*\O_X$ is the Frobenius homomorphism and
that $\sigma_i$ ($i>0$) splits $C$ in the exact sequence
$$
\CD
0 @>>> B^i_X @>>> Z^i_X @>C>> \Omega^i_X @>>> 0
\endCD
$$
The natural perfect pairing $\Omega_X^i\otimes \Omega^{n-i}_X\rightarrow
\Omega^n_X$ gives an isomorphism between
${\cal Hom}_X(\Omega^{n-i}_X, \Omega_X^n)$ and $\Omega^i_X$. It is easy
to check that the homomorphism
$$
F_*\Omega^i_X\rightarrow {\cal Hom}_X(\Omega^{n-i}_X, \Omega_X^n)
\cong \Omega^i_X
$$
given by $\omega\mapsto\varphi(\omega)$, where $\varphi(\omega)(z)=
C(\sigma_{n-i}(z)\wedge \omega)$, splits $\sigma_i$.
\end{pf}

\subsection{Bott vanishing}
Let $X$ be a normal variety
and let $j$ denote the inclusion of the smooth locus $U\subseteq X$.
If the Frobenius morphism lifts to
$W_2(k)$ on $X$, then the Frobenius morphism on $U$ also lifts
to $W_2(k)$.
Define the Zariski sheaf $\tilde{\Omega}^i_X$ of $i$-forms on $X$ as
$j_*\Omega^i_U$. Since $\codim(X-U)\geq 2$ it follows (\cite{Loc},
Proposition 5.10) that
$\tilde{\Omega}^i_X$ is a coherent sheaf on $X$.

\begin{thm}
Let $X$ be a projective normal variety such that $F$ lifts to
$W_2(k)$. Then
$$
\H^s(X, \tilde{\Omega}^r_X\otimes L)=0
$$
for $s>0$ and $L$ an ample line bundle.
\end{thm}
\begin{pf}
Let $U$ be the smooth locus of $X$ and let $j$ denote the inclusion of
$U$ into $X$. On $U$ we have by Theorem \ref{split} a split sequence
$$
0\rightarrow \Omega^r_U\rightarrow F_*\Omega^r_U
$$ which pushes down to the split sequence ($F$ commutes with $j$)
$$
0\rightarrow \tilde{\Omega}_X^r\rightarrow F_*\tilde{\Omega}_X^r
$$
Now tensoring with $L$ and using the projection formula we get
injections for $s>0$
$$
\H^s(X, \tilde{\Omega}^r_X\otimes L)\hookrightarrow
\H^s(X, \tilde{\Omega}^r_X\otimes L^p)
$$
Iterating these injections and using that the Zariski sheaves
are coherent one gets the desired vanishing theorem by
Serre's theorem.
\end{pf}


\subsection{Degeneration of the Hodge to de Rham spectral sequence}

Let $X$ be a projective normal variety with smooth locus $U$.
Associated with the complex
$\tilde{\Omega}^\bullet_X$ there is a spectral sequence
$$
E_1^{pq}=\H^q(X, \tilde{\Omega}^p_X)\implies \H^{p+q}(X,
\tilde{\Omega}^\bullet_X)
$$
where $\H^\bullet(X, \tilde{\Omega}^\bullet_X)$ denotes the hypercohomology
of the complex $\tilde{\Omega}^\bullet_X$. This is the Hodge to de Rham
spectral sequence for Zariski sheaves.

\begin{thm}
If the Frobenius morphism on $X$ lifts to $W_2(k)$, then the spectral
sequence degenerates
at the $E_1$-term.
\end{thm}
\begin{pf}
As complexes of sheaves of abelian groups $\tilde{\Omega}^\bullet$ and
$F_*\tilde{\Omega}^\bullet$ are the same so their hypercohomology
agree. Applying hypercohomology to the split injection (Theorem \ref{split})
$$
\sigma:\bigoplus_{0\leq i}\tilde{\Omega}^i_{X/k}[-i]\rightarrow
F_* \tilde{\Omega}^\bullet_X
$$
we get
\begin{eqnarray*}
\sum_{p+q=n} \dim_k E_\infty^{pq}=
\dim_k \H^n(X, \tilde{\Omega}^\bullet_X)&=&
\dim_k \H^n(X, F_*\tilde{\Omega}^\bullet_X)\geq\\
\sum_{p+q=n} \dim_k \H^q(X, \tilde{\Omega}^p_X)&=&
\sum_{p+q=n} \dim_k E_1^{pq}
\end{eqnarray*}
Since $E_\infty^{pq}$ is a subquotient of $E_1^{pq}$, it follows that
$E_\infty^{pq}\cong E_1^{pq}$ so that the spectral sequence degenerates
at $E_1$.
\end{pf}

\section{Toric varieties}

In this section we briefly sketch the definition of toric varieties
following Fulton \cite{Fulton} and demonstrate how the results
of Section~\ref{flift} may be applied.

\subsection{Convex geometry}
Let $N$ be a lattice, $M =
\Hom_{\Z}(N, \Z)$ the dual lattice, and let $V$ be the real vector
space $V = N \otimes_{\Z} \R$.  It is natural to identify the dual
space of $V$ with $M \otimes_{\Z} \R$, and we think of $N \subset V$
and $M \subset V^*$ as the subsets of integer points.

By a cone in $N$ we will mean a subset $\sigma \subset V$ taking the
form $\sigma = \{r_1 v_1 + \dots + r_s v_s ~ | ~ r_i \geq 0 \}$ for
some $v_i \in N$.  The vectors $v_1, \dots, v_s$ are called generators
of $\sigma$.  We define the dual cone to be $\sigma^{\vee} = \{ u \in
V^* | \forall v \in \sigma: \left< u,v \right> \geq 0 \}$.  One may
show that $\sigma^{\vee}$ is a cone in $M$.  A face of $\sigma$ is any
set $\sigma \cap u^{\perp}$ for some $u \in \sigma^{\vee}$.  Any face
of $\sigma$ is clearly a cone in $N$, generated by the $v_i$ for which
$\left< u, v_i \right> = 0$.

Now let $\sigma$ be a strongly convex cone in $N$, this means that
$\{0\}$ is a face of $\sigma$ or equivalently that no nontrivial
subspace of $V$ is contained in $\sigma$.  We define $S_{\sigma}$ to
be the semi group $\sigma^{\vee} \cap M$.  Since $\sigma^{\vee}$ is a
cone in $M$, $S_{\sigma}$ is finitely generated.

\subsection{Affine toric varieties}
If $k$ is any
commutative ring the semigroup ring
$k[S_{\sigma}]$ is a finitely generated commutative $k$-algebra, and
$U_{\sigma} = \Spec k[S_{\sigma}]$ is an affine scheme of finite type
over $k$.  Schemes of this form are called affine toric schemes.

\subsection{Glueing affine toric varieties}
Let $\tau = \sigma \cap u^{\perp}$ be a face of $\sigma$.  One may
assume that $u \in S_{\sigma}$.  Then it follows that $S_{\tau} =
S_{\sigma} + \Z_{\geq 0} \cdot (-u)$, so that $k[S_{\tau}] =
k[S_{\sigma}]_{u}$.  In this way $U_{\tau}$ becomes a principal open
subscheme of $U_{\sigma}$.
This may be used to glue affine toric schemes together.  We define a
fan in $N$ to be a nonempty set $\Delta$ of strongly convex cones in
$N$ satisfying that the faces of any cone in $\Delta$ are also in
$\Delta$ and the intersection of two cones in $\Delta$ is a face of
each.  The affine varieties arising from cones in $\Delta$ may be
glued together to form a scheme $X_k(\Delta)$ as follows.  If $\sigma,
\tau \in \Delta$, then $\sigma \cap \tau \in \Delta$ is a face of both
$\tau$ and $\sigma$, so $U_{\sigma \cap \tau}$ is isomorphic to open
subsets $U_{\sigma\tau}$ in $U_{\sigma}$ and $U_{\tau\sigma}$ in
$U_{\tau}$.  Take the transition morphism $\phi_{\sigma\tau} :
U_{\sigma\tau} \rightarrow U_{\tau\sigma}$ to be the one going through
$U_{\sigma \cap \tau}$.  A scheme $X_k(\Delta)$ arising from a fan
$\Delta$ in some lattice is called a toric scheme.

\subsection{Liftings of the Frobenius morphism on toric varieties}

Let $X = X_k(\Delta)$ be a toric scheme over the commutative ring $k$
of characteristic $p > 0$.  We are going to construct explicitly a
lifting of the absolute Frobenius morphism on $X$ to $W = W_2(k)$.
Define $X^{(2)}$ to be $X_W(\Delta)$.  Since all the
rings $W[S_{\sigma}]$ are free $W$-modules, this is clearly a flat
scheme over $W_2(k)$.  Moreover, the identities $W[S_{\sigma}]
\otimes_W k \cong k[S_{\sigma}]$ immediately give an isomorphism
$X^{(2)} \times_{\Spec W} \Spec k \cong X$.

For $\sigma \in \Delta$, let $F_{\sigma}^{(2)} : W[S_{\sigma}]
\rightarrow W[S_{\sigma}]$ be the ring homomorphism extending $F^{(2)}
: W \rightarrow W$ and mapping $u \in S_{\sigma}$ to $u^p$.  It is
easy to see that these maps are compatible with the transition
morphisms, so we may take $F^{(2)} : X^{(2)} \rightarrow X^{(2)}$ to
be the morphism which is defined by $F_{\sigma}^{(2)}$ locally on
$\Spec W[S_{\sigma}]$.  This gives the lift of $F$ to $W_2(k)$ and
completes the construction.

\subsection{Bott vanishing and the Danilov spectral sequence}

Since toric varieties are normal we get the following corollary
of Section \ref{flift}:

\begin{thm}
Let $X$ be a projective toric variety over $k$. Then
$$
\H^q(X, \tilde{\Omega}^p_X\otimes L)=0
$$
where $q>0$ and $L$ is an ample line bundle. Furthermore the
Danilov spectral sequence
$$
E_1^{pq}=\H^q(X, \tilde{\Omega}^p_X)\implies \H^{p+q}(X,
\tilde{\Omega}^\bullet_X)
$$
degenerates at the $E_1$-term.
\end{thm}
\begin{remark}
One may use the above to prove similar results in characteristic zero. The
key issue is that we have proved that Bott vanishing and degeneration of
the Danilov spectral sequence holds in any prime characteristic.
\end{remark}


\section{Flag varieties}
In this section we generalize Paranjape and Srinivas result on
non-lifting of Frobenius on flag varieties not isomorphic to $\P^n$. The key
issue is that one can reduce to flag varieties with rank $1$
Picard group. In many of these cases one can exhibit ample line bundles
with Bott non-vanishing.

We now set up notation.
Let $G$ be a semisimple algebraic group over $k$ and fix a Borel subgroup
$B$ in $G$. Recall that (reduced) parabolic subgroups $P\supseteq B$ are
given by subsets of the simple root subgroups of $B$. These correspond
bijectively to subsets of nodes in the Dynkin diagram associated with
$G$. A parabolic subgroup $Q$ is contained in
$P$ if and only if the simple root subgroups in $Q$ is a subset
of the simple root subgroups in $P$. A maximal parabolic subgroup is
the maximal parabolic subgroup not containing a specific simple
root subgroup.

We shall need the following result from the appendix to
\cite{MeSri}

\begin{prop}
\label{splitimplieslift}
If the sequence
$$
0 @>>> B^1_X @>>> Z^1_X @>C>> \Omega^1_X @>>> 0
$$
splits, then the Frobenius morphism on $X$ lifts to $W_2(k)$.
\end{prop}

We also need the following fact derived from
(\cite{Hartshorne}, Proposition II.8.12 and Exercise II.5.16(d))

\begin{prop}
\label{diffilt}
Let $f:X\rightarrow Y$ be a smooth morphism between smooth varieties
$X$ and $Y$. Then for every $n\in \N$ there is a filtration
$F^0\supseteq F^1 \supseteq \dots$ of $\Omega^n_X$ such that
$$
F^i/F^{i+1}\cong f^*\Omega_Y^i\otimes\Omega_{X/Y}^{n-i}
$$
\end{prop}

\begin{lemma}
\label{fibrlemma}
Let $f:X\rightarrow Y$ be a surjective, smooth and projective morphism between
smooth varieties $X$ and $Y$ such that the
fibers have no non-zero global $n$-forms, where $n>0$. Then there is
a canonical isomorphism
$$
\Omega_Y^\bullet\rightarrow f_*\Omega_X^\bullet
$$
and  a splitting
$\sigma:\Omega^1_X\rightarrow Z^1_X$ of
the Cartier operator $C: Z^1_X\rightarrow \Omega^1_X$ induces
a splitting $f_*\sigma:\Omega^1_Y\rightarrow Z^1_Y$ of
$C:Z^1_Y\rightarrow \Omega^1_Y$.
\end{lemma}
\begin{pf}
Notice first that $\O_Y\rightarrow f_*\O_X$ is an isomorphism of rings as
$f$ is projective and smooth. The assumption on the fibers translates into
$f_*\Omega_{X/Y}^n\otimes k(y)\cong \H^0(X_y, \Omega^n_{X_y})=0$ for geometric
points $y\in Y$,
when $n>0$. So we get $f_*\Omega_{X/Y}^n=0$ for $n>0$.
By Proposition \ref{diffilt} this means that all of the natural
homomorphisms $\Omega_Y^n\rightarrow f_*\Omega^n_X$ induced by
$\O_Y\rightarrow f_*\O_X\rightarrow f_*\Omega^1_X$ are isomorphisms giving
an isomorphism of complexes
$$
\CD
0    @>>> \O_Y    @>>> \Omega^1_Y    @>>> \Omega^2_Y    @>>> \dots \\
@.         @VVV           @VVV              @VVV       \\
0    @>>> f_*\O_X @>>> f_*\Omega^1_X @>>> f_*\Omega^2_X @>>> \dots
\endCD
$$
This means that the middle arrow in the commutative diagram
$$
\CD
0 @>>> B^1_Y @>>> Z^1_Y @>C>> \Omega_Y @>>> 0\\
@.       @VVV       @VVV         @VVV  \\
0 @>>> f_* B^1_X @>>> f_* Z^1_X @>f_* C>> f_*\Omega_X @>>> 0\\
\endCD
$$
is an isomorphism and the result follows.
\end{pf}

\begin{cor}
\label{maxparnonlift}
Let $Q\subseteq P$ be two parabolic subgroups of $G$. If the
Frobenius morphism on $G/Q$ lifts to $W_2(k)$, then the
Frobenius morphism on $G/P$ lifts to $W_2(k)$.
\end{cor}
\begin{pf}
It is well known that $G/Q\rightarrow G/P$ is a smooth projective
fibration, where the fibers are isomorphic to $Z=P/Q$. Since $Z$
is a rational projective smooth variety it follows from
(\cite{Hartshorne},  Exercise II.8.8) that $\H^0(Z, \Omega^n_Z)=0$ for
$n>0$. Now the result follows from Lemma \ref{fibrlemma} and
Proposition \ref{splitimplieslift}.
\end{pf}

In specific cases one can prove using the ``standard'' exact sequences
that certain flag varieties do not have Bott vanishing. We go on to do this
next.

Let $Y$
be a smooth divisor in a smooth variety $X$. Suppose that $Y$ is
defined by the sheaf of ideals $I\subseteq \O_X$. Then (\cite{Hartshorne},
Proposition II.8.17(2) and Exercise II.5.16(d)) gives for
$n\in \N$ an exact sequence of $\O_Y$-modules
$$
0\rightarrow \Omega^{n-1}_Y\otimes I/I^2\rightarrow
\Omega^n_X\otimes\O_Y\rightarrow \Omega^n_Y\rightarrow 0
$$
 From this exact sequence and induction on $n$ it follows that
$\H^0(\P^n, \Omega^j_{\P^n}\otimes\O(m))=0$, when $m\leq j$ and $j>0$.

\subsection{Quadric hypersurfaces in $\P^n$}
\label{quadric}
Let $Y$ be a smooth quadric hypersurface in $\P^n$, where $n\geq 4$. There
is an exact sequence
$$
0\rightarrow \O_Y(1-n)\rightarrow \Omega^1_{\P^n}\otimes \O(3-n)\otimes\O_Y
\rightarrow \Omega^1_Y\otimes \O_Y(3-n)\rightarrow 0
$$
 From this it is easy to deduce that
$$
\H^{n-2}(Y, \Omega^1_Y\otimes \O_Y(3-n))\cong
\H^1(Y, \Omega^{n-2}_Y\otimes \O_Y(n-3))\cong k
$$
using that $\H^0(\P^n, \Omega^j_{\P^n}\otimes\O(m))=0$, when $m\leq j$ and
$j>0$.
\subsection{The incidence variety in $\P^n\times \P^n$}
\label{inc}
Let $X$ be the incidence variety of lines and hyperplanes in $\P^n\times
\P^n$, where $n\geq 2$. Recall that $X$ is the zero set of
$x_0 y_0+\dots+x_n y_n$, so that
there is an exact sequence
$$
0\rightarrow \O(-1)\times \O(-1)\rightarrow \O_{\P^n}\times\O_{\P^n}
\rightarrow \O_X\rightarrow 0
$$
Using K\"unneth it is easy to deduce that
$$
\H^{2n-2}(X, \Omega^1_X\otimes \O(1-n)\times \O(1-n))\cong
\H^1(X, \Omega^{2n-2}\otimes \O(n-1)\times \O(n-1))\cong k
$$

\subsection{Bott non-vanishing for flag varieties}

In this section we search for specific maximal parabolic subgroups $P$ and
ample line bundles $L$ on $Y=G/P$, such that
$$
\H^i(Y, \Omega_Y^j\otimes L)\neq 0
$$
where $i>0$.
These are instances of Bott non-vanishing. This will be used in Section
\ref{nonlift} to prove non-lifting of Frobenius for a large
class of flag varieties.

Let $\O(1)$ be the ample generator of $\Pic Y$.
By flat base change one may produce examples of Bott non-vanishing
for $Y$ for fields of arbitrary prime characteristic by restricting
to the field of the complex numbers.
This has been done in the setting
of Hermitian symmetric spaces, where the cohomology groups
$\H^p(Y, \Omega^q\otimes\O(n))$ have been thoroughly investigated by
Sato \cite{Sato} and Snow \cite{Snow1}\cite{Snow2}. We now show that
these examples exist. In each of the following subsections $Y$ will
denote $G/P$, where $P$ is the maximal parabolic subgroup not
containing the root subgroup corresponding to the marked simple
root in the Dynkin diagram. These flag manifolds are the irreducible
Hermitian symmetric spaces.

\subsubsection{Type $A$}
\label{A}
\begin{picture}(30,30)(-20,30)
\put(  0,0){\circle*{4}}
\put(-4.5,-3){$\times$}
\put(  1.5,0){\line(1,0){22}}
\put( 25,0){\circle*{4}}
\put( 30.6,-0.5){.\,.\,.}
\put( 50,0){\circle*{4}}
\put( 51.5,0){\line(1,0){22}}
\put( 75,0){\circle*{4}}
\put( 76.5,0){\line(1,0){22}}
\put(100,0){\circle*{4}}
\put(105.6,-0.5){.\,.\,.}
\put(125,0){\circle*{4}}
\put(126.5,0){\line(1,0){22}}
\put(150,0){\circle*{4}}
\put(145,-3){$\times$}
\end{picture}
\vskip 2.0truecm

If $Y$ is a Grassmann variety not isomorphic to projective
space ($Y=G/P$, where $P$ corresponds to leaving out a simple
root which is not the left or right most one), one
may prove (\cite{Snow1}, Theorem 3.3) that
$$
\H^1(Y, \Omega_Y^3\otimes\O(2))\neq 0
$$

\subsubsection{Type $B$}

\begin{picture}(30,30)(-20,30)
\put(  0,0){\circle*{4}}
\put(  1.5,0){\line(1,0){22}}
\put( 25,0){\circle{3}}
\put( 26.5,0){\line(1,0){22}}
\put( 50,0){\circle{3}}
\put( 55.6,-0.5){.\,.\,.}
\put( 75,0){\circle{3}}
\put( 76.5,0){\line(1,0){22}}
\put(100,0){\circle{3}}
\put(100, 1.5){\line(1,0){25}}
\put(100,-1.5){\line(1,0){25}}
\put(108.2,-3){$>$}
\put(125,0){\circle{3}}
\end{picture}
\vskip 2.0truecm
Here $Y$ is a smooth quadric hypersurface in $\P^n$, where
$n\geq 4$ and Bott non-vanishing follows from Section \ref{quadric}.

\subsubsection{Type $C$}

\begin{picture}(30,30)(-20,30)
\put(  0,0){\circle{3}}
\put(  1.5,0){\line(1,0){22}}
\put( 25,0){\circle{3}}
\put( 26.5,0){\line(1,0){22}}
\put( 50,0){\circle{3}}
\put( 55.6,-0.5){.\,.\,.}
\put( 75,0){\circle{3}}
\put( 76.5,0){\line(1,0){22}}
\put(100,0){\circle{3}}
\put(100, 1.5){\line(1,0){25}}
\put(100,-1.5){\line(1,0){25}}
\put(108.2,-3){$<$}
\put(125,0){\circle*{4}}
\end{picture}
\vskip 2.0truecm
By (\cite{Snow2}, Theorem 2.2) it follows that
$$
\H^1(Y, \Omega_Y^2\otimes\O(1))\neq 0
$$

\subsubsection{Type $D$}

\begin{picture}(30,30)(-20,30)
\put(  0,0){\circle*{4}}
\put(  1.5,0){\line(1,0){22}}
\put( 25,0){\circle{3}}
\put( 26.5,0){\line(1,0){22}}
\put( 50,0){\circle{3}}
\put( 55.6,-0.5){.\,.\,.}
\put( 75,0){\circle{3}}
\put( 76.5,0){\line(1,0){22}}
\put(100,0){\circle{3}}
\put(101.5, 0.5){\line(2, 1){22}}
\put(125,12){\circle*{4}}
\put(101.5,-0.5){\line(2,-1){22}}
\put(125,-12){\circle*{4}}
\end{picture}
\vskip 2.0truecm

For the maximal parabolic $P$ corresponding to the leftmost marked
simple root, Y=$G/P$
is a smooth quadric hypersurface in $\P^n$, where
$n\geq 4$ and Bott non-vanishing follows from Section \ref{quadric}.
For the maximal parabolic subgroup corresponding to one
of the two rightmost marked simple roots we get by
(\cite{Snow2}, Theorem 3.2) that
$$
\H^2(Y, \Omega^4_Y\otimes\O(2))\neq 0
$$

\subsubsection{Type $E_6$}

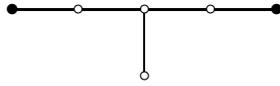
\begin{picture}(30,30)(-20,20)
\put(  0,0){\circle*{4}}
\put(  1.5,0){\line(1,0){22}}
\put( 25,0){\circle{3}}
\put( 26.5,0){\line(1,0){22}}
\put( 50,0){\circle{3}}
\put( 50,-1.5){\line(0,-1){22}}
\put( 50,-25){\circle{3}}
\put( 51.5,0){\line(1,0){22}}
\put( 75,0){\circle{3}}
\put( 76.5,0){\line(1,0){22}}
\put(100,0){\circle*{4}}
\end{picture}
\vskip2.0truecm

By (\cite{Snow2}, Table 4.4) it follows that
$$
\H^3(Y, \Omega^5\otimes\O(2))\neq 0
$$

\subsubsection{Type $E_7$}

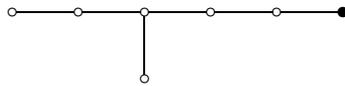
\begin{picture}(30,30)(-20,20)
\put(  0,0){\circle{3}}
\put(  1.5,0){\line(1,0){22}}
\put( 25,0){\circle{3}}
\put( 26.5,0){\line(1,0){22}}
\put( 50,0){\circle{3}}
\put( 50,-1.5){\line(0,-1){22}}
\put( 50,-25){\circle{3}}
\put( 51.5,0){\line(1,0){22}}
\put( 75,0){\circle{3}}
\put( 76.5,0){\line(1,0){22}}
\put(100,0){\circle{3}}
\put(101.5,0){\line(1,0){22}}
\put(125,0){\circle*{4}}
\end{picture}
\vskip 2.0truecm

By (\cite{Snow2}, Table 4.5) it follows that
$$
\H^4(Y, \Omega^6\otimes\O(2))\neq 0
$$

\subsubsection{Type $G_2$}
\label{G}
\begin{picture}(30,30)(-20,30)
\put( 0,0){\circle*{4}}
\put( 0, 1.5){\line(1,0){25}}
\put(  1.5,0){\line(1,0){22}}
\put( 0,-1.5){\line(1,0){25}}
\put(8.2,-3){$<$}
\put(25,0){\circle{3}}
\end{picture}
\vskip2.0truecm
Here $Y$ is a smooth quadric hypersurface in $\P^6$ and Bott non-vanishing
follows from Section \ref{quadric}.
\subsection{Non-lifting of Frobenius for flag varieties}
\label{nonlift}
We now get the following
\begin{thm}
Let $Q$ be a parabolic subgroup contained in a maximal
parabolic subgroup $P$ in the list \ref{A} - \ref{G}. Then
the Frobenius morphism on $G/Q$ does not lift to
$W_2(k)$. Furthermore if $G$ is of type $A$, then the Frobenius
morphism on any flag variety $G/Q\not\cong \P^m$ does not
lift to $W_2(k)$.
\end{thm}
\begin{pf}
If $P$ is a maximal parabolic subgroup in the list \ref{A}-\ref{G}, then
the Frobenius morphism on $G/P$ does not lift to $W_2(k)$. By Corollary
\ref{maxparnonlift} we get that the Frobenius morphism on $G/Q$ does not
lift to $W_2(k)$. In type $A$ the only flag variety not admitting
a fibration to a Grassmann variety $\not\cong\P^m$ is the incidence
variety. Non-lifting of Frobenius in this case follows from Section
\ref{inc}.
\end{pf}

\begin{remark}
The above case by case proof can be generalized to include
projective homogeneous $G$-spaces with non-reduced stabilizers.
It would be nice to prove in general that the only flag variety
enjoying the Bott vanishing property is $\P^n$. We know of no other
visible obstruction to lifting Frobenius to $W_2(k)$ for flag varieties
than the non-vanishing Bott cohomology groups.

\end{remark}

\newpage
\bibliographystyle{amsplain}
\ifx\undefined\bysame
\newcommand{\bysame}{\leavevmode\hbox to3em{\hrulefill}\,}
\fi

\end{document}